\documentclass[10pt,a4paper]{amsart}
\usepackage{todonotes}
\newtheorem{prop}{Proposition}[section]

\newtheorem{theorem}[prop]{Theorem}
\newtheorem{Lemma}[prop]{Lemma}
\newtheorem{Conj}[prop]{Conjecture}
\newtheorem{Def}[prop]{Definition}

\newcommand\bH{\mathbb{H}}
\title[congruence subgroups]{Optimal Farey sequence  for the Congruence subgroup $\Gamma_0(2^{n})$
}
\numberwithin{equation}{section}

 \author {Nhat Minh Doan}
 \address{Department of Mathematics, National University of Singapore, Singapore \& Institute of Mathematics, Vietnam Academy of Science and Technology, Vietnam}
 \email{minh.dn@nus.edu.sg, dnminh@math.ac.vn}
 
\author{Sang-hyun Kim}
 \address{School of Mathematics, Korea Institute for Advanced Study (KIAS), Korea}
 \email{kimsh@kias.re.kr}
 
\author{Mong Lung Lang}
 \address{Singapore}
 \email{lang2to46@gmail.com}
 
\author{Ser Peow Tan}
 \address{Department of Mathematics, National University of Singapore, Singapore}
 \email{mattansp@nus.edu.sg}
\begin{document}

\baselineskip12pt

\maketitle
\vspace{-0.3in}

\begin{abstract}
  We prove that
  $\Gamma_0(2^n)$ ($n\ge2$) has  a Farey sequence $\{e_i\}$ such that  $e_i \le
   2^{n-1}
  $
    for all $e_i$.
     The above upper bound is optimal
      and there exists a unique $j$ such that  $e_j= 2^{n-1} $. For each $e_i$, there exists a 
       unique $a_i$ such that $\{ a_i/e_i\}\cup \{\infty\}$ is the set of ideal vertices of a fundamental domain
    of $\Gamma_0(2^n)$ whose side-pairings give a set of independent generators of $\Gamma_0(2^n)$.

\end{abstract}

\section{Introduction}

\subsection{ The main result
}\label{ss:Mainresult}
Let  \begin{equation}\label{eqn:sequenceT} T= \{0/1=a_1/e_1,\ldots\ldots, a_m/e_m = 1/1\} \cup \{\infty\}
\end{equation}
be the ordered set of cusps (ideal vertices)   of a special polygon $P$ (a special type of fundamental domain) of $\Gamma_0(2^n)$ ($n\ge2$). The side-pairings of $P$ give an independent set of generators for $\Gamma_0(2^n)$. We call the denominator sequence $\{e_i\}$ a {\it Farey sequence} for $\Gamma_0(2^n)$. 
Motivated by the work of Kulkarni \cite{Ku},
we believe that  $\{e_i\}$ sheds light on some interesting and important geometric and arithmetic properties of $\Gamma_0(2^n)$, if $P$ is constructed carefully. In  this  note, we prove:
 \begin{theorem}\label{thm:main}(Theorem 8.1)
     $P$ can be
  constructed in such a way that $\{e_i\} $ satisfies
  \begin{equation}\label{eqn:boundonFareyseq}
      e_i\le 2^{n-1}\,\mbox{ for all $i$ and $e_j = 2^{n-1}$ for a unique $j$}.
  \end{equation}
  \end{theorem}
  For each $e_i$, there exists a 
 unique $a_i$ such that $\{ a_i/e_i\}\cup \{\infty\}$ is the set of cusps of a special polygon 
 of $\Gamma_0(2^n)$.
 
 This is achieved in subsection 8.1. Note that the upper bound $2^{n-1}$ is optimal, which confirms our conjecture for $\Gamma_0(N)$ in this special case; see \cite[Conjecture 3.16]{DKLT}.
 \subsection{A weaker result}\label{ss:weakerresult}    A large part of our method works for $\Gamma_0(p^n)$ when $p$ is an odd prime and $n\ge 2$. However, in the last part of the construction, we were unable to obtain an optimal upper bound for $\{e_i\}$.  The best we can get is 
  the following.
 
 \smallskip
 \begin{theorem}\label{thm:weakerthm}(Theorem 10.1)
      { Let $p>2$ be  prime and $n\ge2$. Then
 	$\Gamma_0(p^n)$ has  a Farey sequence $\{e_i\}$ such that 
 	\begin{equation}\label{eqn:weakerbound}
 	    e_i\le(f_{p+1}/2) p^{n-1} 
 	+ f_pp^{n/2}\,\,for\, all\, i,
 	\end{equation}
 where $f_m$ is the $m$-th Fibonacci number} ($f_1=f_2=1, f_m=f_{m-1}+f_{m-2}$ for all $m \in \mathbb Z$).
 \end{theorem}

 \subsection{Previous results}\label{ss:prevresults} The idea of constructing a special type of fundamental domain $P$ for $\Gamma_0(N)$ ({\it special polygon}) where the sequence of cusps of $P$ are Farey neighbors and the side-pairings give an independent set of generators for $\Gamma_0(N)$ goes back to Kulkarni in \cite{Ku}. This was developed further in \cite{CLLT} where an algorithm for the construction was given, eliminating the trial and error ambiguity in the original construction. This has been implemented by others for example in \cite{Mag} and \cite{Sag}. However it does not necessarily give optimal bounds for the denominators. Special polygons with optimal bounds for the Farey sequence $\{e_i\}$ for the cases where $N=p$ or $ p^2$, $p$  prime  were constructed in \cite{DKLT}. 
 
 \subsection{Farey sequences}\label{ss:Fareyseq}  In \cite{CLLT}, the (ordered) collection of members in  $T$ in (\ref{eqn:sequenceT})   was   referred to as  a 
 Farey sequence of $\Gamma_0(2^n)
 $. Since the $a_i$'s are uniquely determined by the $e_i$'s (see subsection \ref{ss:FCforfractions}), and we are mostly interested in the 
 properties of the $e_i$'s,  there is no point in constructing the sequence $\{a_i/e_i\}$.
  We therefore only construct the sequence $ \{e_i\}$.  In this paper, we call this a Farey sequence of $\Gamma_0(2^n)$, and more generally, the set of denominators of a sequence of Farey neighbors between $0$ and $1$ is called a Farey sequence. Similarly, by abuse of notation, we also denote by $F_r$ the ordered sequence of denominators of the Farey sequence of order $r$.  An obvious advantage is that expressions such as 
    (5.1) and (5.5) are neater and  would look messy in the old notation (each entry would take the form  $(s+vt)/(a+vb_0p^m)$ instead of  $a+vb_0p^m$).

  \subsection{Notations}\label{ss:Notations} 
  Let $\Gamma:=PSL(2,\mathbb Z)$ which acts by M\"obius transformations isometrically on the hyperbolic upper half plane $\mathbb H$.  Translates of the ideal triangle with vertex set $\{\infty, 0,1\}$ by $\Gamma$ are called {\it Farey triangles} and translates of the infinite line $(\infty, 0)$ by $\Gamma$ are called even lines. The ideal endpoints of an even line are {\it Farey neighbors}. Translates of the triangle with vertex set $\{\infty, 0, e^{i\pi/3} \}$ by $\Gamma$ are called {\it special triangles}. Note that each Farey triangle contains $3$ special triangles. $\Gamma_0(N)$ is the subgroup of $\Gamma$ consisting of all matrices whose $(2,1)$ entry is divisible by $N$. Denote by $\Sigma(N)$ the surface $\mathbb{H}/\Gamma_0(N)$. Note that $\Sigma(N)$ is tessellated by Farey triangles and special triangles (if $\Sigma(N)$ has order $3$ orbifold points). The natural collar of a cusp $x$ of $\Sigma(N)$, denoted by $C(x)$ is the collar  about $x$ bounded by the horocycle of length $2$. Note that the natural collars of distinct cusps are disjoint.  The projection map of $\mathbb{H}$ to $\Sigma(N)$ is denoted by $\pi$.
  
  Let $T$ be given as in (1.1). A pair of adjacent members is called a side. Adjacent terms will be Farey neighbors by construction. In our new notation, $(e_i, e_{i+1})$ is a 
   side and  $e_i$ is a cusp (vertex). However,  when the $e_i$'s are replaced by numbers, one has to be careful.
  Take  the special polygon $\{a_i/e_i\}\cup\{\infty\}= \{0/1,1/3,1/2,2/3,1/1\}\cup\{\infty\} $ of $\Gamma_0(6)$ for instance.
  Following our convention (subsection \ref{ss:Fareyseq}),   
   $ F=  \{1,3,2,3,1\} $ is  a Farey sequence of  $\Gamma_0(6)$. The pairs
    $(1,3)$, $(3,2)$,  $(2,3)$ and  $(3,1)$ are sides. We will also use the convention $(1,3) \in F$ to denote that $(1,3)$ is a side of $F$, this should not cause any confusion.
   It is not correct to say $3$ is a cusp
    (vertex) since it is not unique. The
    correct way is to say the cusp $3$ in $(1,3)$ and 
     the cusp 3 in  (3,1).  These are two different cusps. Where the need arises, we will still use $x/y$ to denote a cusp and $[x/y]$ to denote its equivalence class.

     \medskip
     
      Let $(a,b)$ be a side of a Farey sequence $\{e_i\}$. Borrowing from the Farey construction, 
      \begin{equation}\label{eqn:Fareyaddition}
          a\oplus b  \mbox {  denotes the number   $a+b$ which we insert into  }  (a,b).\end{equation}

   \subsection{Construction of the optimal Farey sequence $\{e_i\}$ for $\Gamma_0(2^n)$.}\label{ss:constrFS}  This is achieved in Section \ref{s:mainconstruction} where
     $\{e_i\}$ is  constructed. 
    In short, $\{e_i\}$ is a union given as follows:
   \begin{equation}\label{eqn:Soneunion}\{e_i\}= S_1 \cup \{a\oplus b\,:\,  a+b=z2^i, \, i> n/2\} \qquad (\mbox{see  \ref{eqn:Fareyaddition}} )
     ,\end{equation}
  \noindent  where $S_1$ is the union of $F_r$, $r=\lfloor 2^{n/2}\rfloor$, and some chains we call $N_e$
   (see (\ref{eqn:optimalchains})).
  In detail,
   $F_r$ is given in Section 2,  the chains are given in (\ref{eqn:optimalchains}) and $T\setminus S_1$ is given in (\ref{eqn:seqforT}). 
    In Section \ref{s:64},   we construct a Farey sequence of $\Gamma_0(64)$ and demonstrate how
     one can work out the general case similarly. It is understood that members in $\{e_i\}$ are ordered
      according to the Farey construction arising from the operation in (\ref{eqn:Fareyaddition}) (see also subsection \ref{ss:FareySeqGamma0N}).
   
   \subsection{Upper bounds  for Farey sequences for $N=2^n$.}\label{ss:upperbounds} 
   An upper bound for  vertices of $F_r$ is $r:=\lfloor 2^{n/2}\rfloor$, essentially by definition
   (see Section \ref{s:ConstFareyseq}).  An upper bound for  vertices of $S_1$ is $2^{n-1}$ (see Lemma \ref{lem:boundsforchains}).
 An upper bound for  vertices of $T\setminus S_1$ is also $2^{n-1}$ (see Theorem \ref{thm:mainthm}).

 \subsection{Geometric discussion}\label{ss:Discussion} We give a  geometric interpretation of our construction here which may be helpful to the more geometrically inclined reader. Recall that we want to construct a special type of fundamental domain for the surface $\Sigma(N):=\mathbb{H}^2/\Gamma_0(N)$. In our earlier study \cite{DKLT} on optimal special polygons for $\Sigma(N)$ where $N=p$ or $p^2$ ($p$ prime), the first important observation was that a special polygon can be chosen to contain the ideal polygon with vertices $F_r$, where $r:=\lfloor \sqrt N\rfloor$.
 This still holds in our case  and is the first step in our construction. The second important fact used was that it was only necessary to add isolated Farey triangles or special triangles to complete the special polygon. This no longer holds here ($N=p^n$, $n>2$). The main difference is that in general, $\Sigma(N)$  has cusps with intermediate widths, (see section \ref{ss:cusps}) whereas, for $N=p$ or $p^2$ there were only cusps with minimal or maximal widths ($1$ or $N$). Recall that the cusps of $\Sigma(N)$ correspond to the $\Gamma_0(N)$-orbits of $x/y \in \mathbb{Q} \cup \{\infty\}$, and that the natural collars of the  cusps  are mutually disjoint. Call two (not necessarily distinct) cusps of $\Sigma(N)$ {\it Farey-adjacent} if they are connected by an even line. We summarize the steps/observations  towards proving our result as follows (see subsection \ref{ss:cusps} for (2),(3), and (4)):
 \begin{enumerate}
    \item We may start with $F_r$ in the construction of the Farey sequence $\{e_i\}$ for $\Sigma(p^n)$;
    
     \item The cusp $[0/1]$ is the unique cusp of maximal width and it is Farey-adjacent to all cusps, including itself;
     
     \item A  cusp $[x/y]$ of minimal or intermediate width is Farey-adjacent only to the cusp $[0/1]$;
     
      \item The collar neighborhood of every cusp of intermediate width is already partially filled in by $F_r$. That is, the projection $\pi$ of the ideal polygon with cusp set $F_r$ to $\Sigma(N)$ intersects the collars of all cusps of intermediate width;
     
     \item There is a systematic way to fill in $C(x)$, the collars of cusps $x$ of intermediate widths  sequentially, by inserting Farey triangles to $F_r$, while bounding the denominators by $p^{n-1}$. This is by appending the chains $N_{\epsilon}$ to $F_r$ to get $S_1$ (see section \ref{s:Sone} for the construction of $S_1$ and section \ref{s:Upperbounds} for the bounds);
     
     \item We complete the construction by filling in the collar  $C([0/1])$. In the process, all other cusps  would also be filled in by (3) above, and we obtain a Farey sequence for $\Gamma_0(p^n)$;
     
     \item In the case $p=2$, in the last step we only need to fill in isolated Farey triangles, whose denominators are still bounded by $2^{n-1}$. If $p \neq 2$, it may be necessary to fill in blocks with at most $2p-2$ Farey triangle pieces. This leads to the bound in Theorem \ref{thm:weakerthm}.
     
 \end{enumerate}

   
\noindent {\it Remark:}  Throughout this  article, when we consider $\Gamma_0(p^n)$, $n\ge 3$  unless otherwise mentioned.\\

\noindent
\textbf{Acknowledgment.} Nhat Minh Doan is funded by the Singapore National Research Foundation (NRF) under grant E-146-00-0029-01. Sang-hyun Kim is supported by Mid-Career Researcher Program (RS-2023-00278510) through the National Research Foundation funded by the government of Korea. Ser Peow Tan is supported by the National University of Singapore academic research grant A-8000989-00-00.

\section {Constructing Farey sequences}\label{s:ConstFareyseq}

\subsection{The Farey sequence $F_r$}\label{ss:FareySeqFr} Let $N>3$ and $r=\lfloor \sqrt{N}\rfloor$. We are mostly interested in the case $N=p^n$ where $p$ is prime (in particular, when $p=2$) and $n \ge 3$.
Let $F_r$ be the ordered sequence of denominators of the Farey sequence  of order $r$ (rationals between $0$ and $1$ with denominators $\le r$).  
Members of $F_r$ are ordered by the {\em Farey construction}
of Farey sequences (see subsection \ref{ss:FareySeqGamma0N}). To be precise,
$F_{r}$ can be  constructed inductively as follows:

\begin{enumerate}
	\item[(i)]
	Set
	$D_1= \{1,1\}.$  The pair $(1,1)$ is called a side of $D_1$.
	
	\vspace{.1cm}
	\item[(ii)]  Insert $1+1=2$ into $D_1$.  This gives $D_2 = \{ 1,2,1\}$.
	Adjacent pairs $(1,2)$ and $(2,1)$ are the sides of $D_2$. Note that  $(1,1)$ is no longer a 
	 side of $D_2$.

		\vspace{.1cm}

  \item[(iii)] For $k \ge2$, let $(a,b)$ be a side of $D_k$ such that $a+b \le r$. Insert  $a+b$ into $(a, b)$. This gives $D_{k+1} = \{ \ldots, a, a+b, b, \ldots\}$.  
	  Note that $(a, b)$ is no longer a side of $D_{k+1}$. For example, suppose that $(a, b)= (1,2)$ in $D_2$.
	   Then $D_3 = \{ 1, 3,2,1\}$.

	  \vspace{.1cm}
	  \item[(iv)]  This process will terminate in the sequence $D_q:=F_r$, when for all sides $(a,b)$ of $D_q$, $a+b >r$.
	  
	\end{enumerate}

\subsection{A Farey sequence for $\Gamma_0(N)$}\label{ss:FareySeqGamma0N} The algorithm in subsection \ref{ss:FareySeqFr} can be regarded as the first step towards constructing a Farey sequence (hence special polygon) for $\Gamma_0(N)$. Indeed, if  we replace (iii) in the algorithm above  by
 \bigskip
\begin{enumerate}
    \item [(iii$'$)]  For $k \ge 2$, let $(a,b)$   be a side of $D_k$ such that  $a^2 + b^2$ and $a^2+ab+b^2$
	 are not multiples of $N$;  $ax+by$ is not a multiple of $N$  for every side $(x, y)$ of $D_k$; and $a+b$ is minimal among the sides $(a,b)$ with these properties.
	  Insert  $a+b$ into $(a, b)$. This gives $D_{k+1} = \{ \ldots, a, a+b, b, \ldots\}$.  
\end{enumerate}

Then the algorithm using (i), (ii) and  (iii$'$) terminates after a finite number of steps in a Farey sequence $\{e_i\}:=D_s$ for $\Gamma_0(N)$, since $\Gamma_0(N)$ has finite index in the modular group $\Gamma$ (see \cite{CLLT}). We call this algorithm $(CLLT)$.

The first key observation in \cite{DKLT} was that the sequence $F_r$  ($r=\lfloor \sqrt{N}\rfloor$) always occur as an intermediate sequence in the algorithm, and so in fact, we can start the algorithm from the sequence $F_r$ instead of $D_1$, that is, replace (i) with (i$'$) Set $D_1=F_r$ (step 1 in section \ref{ss:Discussion}). This comes from the minimality condition of (iii$'$). In fact, this minimality condition is not necessary, and was not assumed in \cite{CLLT}, and one can still get a Farey sequence for $\Gamma_0(N)$ if we drop it.

\smallskip

The special polygon with the side pairings can be constructed from this sequence and the pairing data. We first start with a definition:

\bigskip
\begin{Def}\label{def:defofsides}
The sides $(a,b)$ of $D_s$ satisfying $a^2+b^2$ is a multiple of $N$ are called {\it even sides}, the sides satisfying $a^2+ab+b^2$ is a multiple of $N$ are called {\it odd sides}. Collectively, even and odd sides are called {\it elliptic sides} as the corresponding side pairings are elliptic. Two sides $(a,b)$ and $(x,y)$ of $D_s$ satisfying $ax+by$ is a multiple of $N$ are 
called  {\it paired sides}. For the intermediate sequences $D_r$ where $r<s$, a non-elliptic side of $D_r$ which is not paired is called a {\it free side}.
\end{Def}

\noindent{\bf Remark:} For a given $N$, the definition of an elliptic side makes sense for denominators of any even line, and similarly one can talk about pairing of $(a,b)$ with $(x,y)$ for co-prime integers $a,b$ and $x,y$ which corresponds to the denominators of a pair of even lines.

\smallskip

 When the process terminates at $D_s$, all sides are either elliptic or paired as there are no more free sides to insert in. There is a unique ideal polygon $P$ (see subsection \ref{ss:FCforfractions} and Remark 3.2) whose denominator sequence is given by $D_s$. 

 
	 



	\subsection{ Farey construction for fractions and Farey sequences}\label{ss:FCforfractions} Let $A_1=\{ 0/1, 1/1\}$. $A_1 \cup \{\infty\}$ are the vertices of a Farey triangle. Following the Farey 
	 construction given in \cite{CLLT},  $ A_{k+1} $ is obtained from $A_k$ by inserting
	 the Farey sum (see (\ref{eqn:deffareyaddition})) of $a/b$ and $c/d$, where $e= (a/b, c/d)$ is a free side of $A_k$.
	 Note that a free side of $A_k$ is a non-elliptic side that is not paired with any non-elliptic side of $A_k$.
	The Farey sum of $a/b$ and $c/d$  is defined as follows:
	  \begin{equation}\label{eqn:deffareyaddition}\frac{a}{b}  \oplus  \frac{c}{d}   = \frac{a+c}{b+d}\,.\end{equation}
	  One sees easily  that  $D_1$ in  subsection \ref{ss:FareySeqFr}
	  is the denominator of $A_1$ and that in general,
	  $D_k$ is the set of denominators of  $A_k$. This is best illustrated by the following example.
	  $$A_2 =\{ 0/1, 1/2, 1/1\},\,\, D_2 =\{ 1, 2,1\},\,
	  A_3=\{ 0/1,  1/3, 1/2,1/1\} \mbox{  and } D_3= \{ 1,3,2,1\}.$$
	 Inserting a Farey sum is equivalent to attaching a Farey triangle. Hence, $D_k$ is the ordered set of denominators of the vertices of an ideal polygon consisting of Farey triangles. A sequence $D_k$ obtained inductively in this way, without any restrictions on the insertions, is called a {\it Farey sequence}.  The following holds (see [CLLT] for (ii) and (iii)).

	\begin{enumerate}
	\item[(i)]	$A_k$ determines uniquely a  Farey sequence $D_k$ and vice versa.
	
	\item[(ii)]	 If $e= (c/a, d/b)$ and $f=(w/x, z/y)$ are sides of $A_k$ such that $e$ and $f$ are paired by an element $M$ of $\Gamma$,
	 then the (2,1)-entry of $M$   is a multiple of $N$, so $M \in \Gamma_0(N)$. Note that the 
	 (2,1)-entry of $M$ is  $ax+by$.

		\item[(iii)]	 If $e= (c/a, d/b)$  is elliptic, then $e$ is self-paired by an elliptic element $M$ of $\Gamma_0(N)$
		whose (2,1)-entry is a multiple of $N$. Note that the (2,1)-entry of $M$ is 
		 either  $a^2+b^2$ or $a^2+ab+b^2$. In the former case, the side is regarded as two semi-infinite sides which are paired, and in the latter case, a special triangle is attached to the side with the resulting semi-infinite sides paired.

\end{enumerate}






\section{Farey sequences of $\Gamma_0(p^n)$ and cusps of $\Sigma(p^n)$}

\subsection{A Farey sequence of $\Gamma_0(p^n)$}
The construction in Subsection \ref{ss:FareySeqGamma0N} terminates in a Farey sequence $D_s$ for $\Gamma_0(p^n)$ (when $N=p^n$) with the following properties:
\begin{enumerate}
	\item[(i)]
	$D_s$ has $1+\lfloor[\Gamma : \Gamma_0(p^n)]/3\rfloor$  members (follows from the index  $[\Gamma:\Gamma_0(p^n)]$).
	
	\vspace{.1cm}
	\item[(ii)]
	For each non-elliptic  side $s=(a, b)  $ of $D_s$, there exists a unique side $s^*=(x,y)$
	of $D_s$
	such that $s$ and $s^*$ are paired, that is,  $ax+by$ is a multiple of $p^n$. (see \cite{Ku} and \cite{CLLT}).  
\end{enumerate}

\smallskip

\begin{Def}\label{def:subsequences}
    A Farey sequence $F$ is called a
{\em Farey sequence} of $\Gamma_0(p^n)$ if $F_r \subseteq F$ and $F$ satisfies  (i) and (ii) above.  A Farey sequence $S$ is called a subsequence of $\Gamma_0(p^n)$
if $F_r\subseteq S$,  and $S$ is a subsequence of a Farey sequence $F$ of $\Gamma_0(p^n)$. 
\end{Def}

It follows from the definition that a Farey subsequence $S$  can be extended to a Farey sequence $F$  of $\Gamma_0(p^n)$ by the Farey construction. Also, the intermediate sequences $D_r$ where $r<s$ of the algorithm using (i$'$), (ii), (iii$'$) in subsection \ref{ss:FareySeqGamma0N} are Farey subsequences of $\Gamma_0(p^n)$.

\smallskip
\noindent {\bf Remark 3.2.}   The sequence  $D_s =\{ e_i \}$  determines 
a special polygon of $\Gamma_0(p^n)$ (see subsection 1.3 and subsection \ref{ss:FCforfractions}). 
For each $e_i$, there exists a 
unique $a_i$ such that $\{ a_i/e_i\}$ is the set of cusps of a special polygon (fundamental domain)
of $\Gamma_0(p^n)$. Special triangles are attached to the odd sides of the ideal polygon with the sides meeting at the finite vertex paired. An even side is paired to itself by the $\pi$ rotation about the translate of $\sqrt{-1}$ lying on it. The sides of  paired sides are paired to each other.



\smallskip

\noindent {\bf Remark 3.3.} The algorithm of subsection \ref{ss:FareySeqGamma0N} inserts $a+b$ into a free side $s=(a,b)$ of a Farey subsequence of $\Gamma_0(p^n)$, chosen so that $a+b$ is minimal among all possible free sides, to obtain a new Farey subsequence for $\Gamma_0(p^n)$. This terminates after a finite number of steps in a Farey sequence of $\Gamma_0(p^n)$, since the index $[\Gamma:\Gamma_0(p^n)]<\infty$. If we do not include the condition that $a+b$ is minimal, we will still terminate in a (possibly different) Farey sequence for $\Gamma_0(p^n)$.

\subsection {  The Farey sequences $F_{[\sqrt 8\,]}$ and  $F_{\sqrt {16}}$  and
	Farey sequences of $\Gamma_0(8)$  and   $\Gamma_0(16)$ }

It is easy to see that 
$F_{[\sqrt 8\,]} =\{1,2,1\}$ and $F_{\sqrt {16}}=\{1,4,3,2,3,4,1\}.$
Further, a simple check shows that the following are Farey sequences for $\Gamma_0(8)$  and $\Gamma_0(16)$  respectively:
\begin{equation}\label{eqn:example}\{1,4,3,2,1\},\,\,\{1,4,3,8,5,2,3,4,1\}.  \end{equation}

\subsection{ Cusps   of $\Sigma(p^n)$ }\label{ss:cusps}


  Recall that, in general,
$\Sigma(N):=\bH/\Gamma_0(N)$  has $\sum _{ d|N} \phi(
\mbox{gcd}\,((d, N/d))$ inequivalent cusps. These cusps correspond to the orbits of the action of $\Gamma_0(N)$ on $\mathbb Q \cup\{\infty\}$, so we also denote a cusp of $\Sigma(N)$ by the equivalence class $[x/y]$. An important invariant of a cusp is:

\begin{Def}\label{def:widthofcusps}
    Let $[x/y]$ be a cusp of $\Sigma(N)$ and let $A\in PSL(2,\mathbb Z)$ be a matrix that
sends $\infty $ to $x/y$. The width  of $[x/y]$ is the smallest
positive integer $w := w([x/y])$  such that
\begin{equation}\label{eqn:widthsofcusps}A
\left (\begin{array} {cc}
	1&w\\
	0&1\\
\end{array}
\right ) A^{-1} \in \Gamma_0(N).\end{equation}
\end{Def}

Geometrically,  $w([x/y])$ is the number of even lines intersecting the collar $C([x/y])$. Note that $w([\infty])=1$. We abbreviate $w([x/y])$ by the notation $w(x/y)$.
The following  facts about the  cusps of $\Sigma(N)$ will be useful (see pp. 25 of [Sh]):
\begin{Lemma}\label{lem:cusps1}
Consider $\Gamma_0(N)$ and the associated surface $\Sigma(N)$.
\begin{enumerate} 
    \item[(a)] For each positive divisor $d$ of $N$,
$x/d$ and $y/d$ (in reduced form) are equivalent to each other  if and only if
$x-y$ is a multiple of gcd$\,(d, N/d)$.

    \item[(b)] $x/y$ (in reduced form) is
equivalent to
$x'/gcd\,(y, N)$ for some $x'$. 
\end{enumerate}
\end{Lemma}

Furthermore, if $N=p^n$, we have:
\begin{Lemma}\label{lem:csps2}
Consider $\Gamma_0(p^n)$ and the associated surface $\Sigma(p^n)$, where $p$ is prime. Let $x=a/bp^i$ (in reduced form),  where gcd$\,(ab,p)=1$, and $0 \le i \le n$. 
\begin{enumerate}
	\item [(i)]
	$w(x) = 1$ if $i\ge n/2$  and $w(x) = p^{n-2i}$ if $i\le   n/2$. In particular, $w(x)=p^n$ if and only if $[x] = [0/1]$. Cusps of width $p^{n-2i}$ where $0<i< n/2 $ are called cusps of intermediate widths.
	
	\item[(ii)]  Following the construction of $F_r$ (see subsection \ref{ss:FareySeqFr}) and from Lemma \ref{lem:cusps1}(b),
	 if $i\le n/2$, then
	$x=a/bp^i$ is equivalent to a vertex of $F_r$. 
	\item[(iii)] Every cusp $[x/y] \neq [0/1]$ is Farey-adjacent only to $[0/1]$ and conversely, $[0/1]$ is Farey-adjacent to all cusps.
        \item[(iv)] Suppose that $[x_1]=[1/1], [x_2]=[1/2], \ldots, [x_{p^n}]=[1/p^n]$ are the cusps Farey-adjacent to $[0/1]$ in cyclic order, then  $[x_i]$, has intermediate or minimal width if and only if $p$ divides $i$. Furthermore, if  $x_i$ has width $1$ and $n \ge 3$, then $w(x_{i-p})>1$ and $w(x_{i+p})>1$.
	
	\item[(v)] Suppose that  $N=2^n \ge 4$.
	Let $c = x/y2^{n-1},$ where gcd$\,(xy,2)=1. $
	It follows from the above that  $c$ has width 1 and
	$\Sigma(2^n)$ has a unique cusp of this form.

\end{enumerate}

\end{Lemma}

\section{Farey sequences of $\Gamma_0(64)$  and $\Gamma_0(2^n)$}\label{s:64}
 
 We construct a Farey sequence  for $\Gamma_0(64)$ and demonstrate that  one can 
 construct a Farey sequence for $\Gamma_0(2^n)$ in the same manner (see subsection \ref{ss:generalcase}).
 Hopefully, this  informal  treatment will make our  presentation in Section \ref{s:Sone} and subsection \ref{ss:pequaltwo}
  more friendly and readable.
It is therefore safe to say that the $\Gamma_0(64)$ case is pivotal  (see also Conjecture  \ref{Conj:bound}).
Note first that
\begin{equation}\label{eqn:sequencefor64}F_{\sqrt{64}} =F_8=\{ 1,8,7,6,5,4,7,3,8,5,7,2,   7,5,8,3,7,4,5,6,7,8,1\}.\end{equation}
\subsection{First observation.}\label{ss:observation}   The side $(1,8)$ pairs with $(8,7)$.  In general, if $N=p^{2m}$, the side
 $(a, p^m)$ of $F_r$ pairs with the side $(p^m, p^m-a)$ of $F_r$ ($r=m$).

 \subsection{Constructing chains}\label{ss:chains} We now look at the side $ e_1=(7,6)$. Note that the second entry is a multiple of 2.  One checks that 
 $(7,6)$ is free, and we may insert cusps into $(7,6)$. At each step, we insert a cusp next to $6$ when possible. This gives us  the following chains:
 \begin{equation}\label{eqn:chainsfor64nonoptimal} \{7,13,6\} ,\,\,  P= \{7,13,19,6\},\,\,
N'=  \{7,7+6, 7+2\cdot 6 ,7+3\cdot 6,6\}=\{7,13,19,25,6\} .\end{equation}
The chain stops at $N'$ as $(25, 6)$ pairs with $(6,7)\in F_{\sqrt{64}}$  (here $25\cdot 6+6\cdot7= 192$ is a multiple of 64). The largest member of $N'$ is $25$ which is quite large but still less than $2^5=32$. We can make the largest element smaller by partially filling from the right from the other side too, and doing it in a balanced way. Here we replace the chain $N'$ by
\begin{equation}\label{eqn:chainsfor64optimal}\{7,13,19,6\} \cup \{6,13,7\}, \end{equation}
where $(19,6)$ pairs with $(6,13)$.

\smallskip

 In general, for $N=p^n$, one works on free sides of $F_r$ of the form
$e= (a, bp^m)$ and first constructs  the  chain
\begin{equation}\label{eqn:chainsingeneralnonoptimal} N_e'= \{a,    a+ bp^m,  a+2bp^m,\ldots,  a+vbp^m,         bp^m\} \end{equation}
where $v$ is the smallest positive integer such that  
$(a+vbp^n, bp^m)$ is paired to some side $(x,y)$ of $F_r$. This always exists as the width of the cusp is finite. Indeed, an immediate bound for $v$ is $w-1$ where $w$ is the width of the cusp.
We show in subsection \ref{ss:upperboundsforchains} that $a+vbp^m<p^{n-1}$.

This can then be replaced by the more optimal chain
\begin{equation}\label{eqn:chainsingeneraloptimal}  N_e= \{a,    a+ bp^m,  a+2bp^m,\ldots,  a+v_1bp^m,         bp^m\}\cup \{x,v_2x+y, \ldots, 2x+y, x+y,y \} \end{equation}
where $v_1+v_2=v$, $v_1=v_2$ if $v$ is even and $v_1=v_2+1$ if $v$ is odd. The side $(a+v_1bp^m,  bp^m)$ pairs with $(x,v_2x+y)$.
\smallskip

Geometrically, we are filling in the gaps in the collars of the cusps of intermediate widths by adding Farey triangles to the partial fundamental domain corresponding to $F_r$. Every gap corresponds to  free sides of the form $e= (a, bp^m)$ and $e=(xp^m, y)$ of $F_r$, which are the left and right boundaries of this gap. So after adding in all the chains, we have filled in all the gaps in the collars of the cusps of intermediate widths. Note also that by Lemma \ref{lem:csps2}(iii) the filling in of the chains $N_e$ are independent of each other. We fill in the gaps in a balanced way to reduce the size of the denominators obtained; this will be important in obtaining the optimal bound for $\Gamma_0(2^n)$.





\subsection{The Farey sequences $S_1$ and $T$.}\label{ss:SeqSone}  Let $S_1$ be the union of  $F_{\sqrt {64}}$ and  $N_{e_i}$,
 where $e_i \in\{ e_1= (7,6), e_2=(5,4), e_3=(7,2), e_4=(7,4), e_5=(5,6)\}$  (the free sides $(a,b)$ of $F_{8}$, where $b$ is even). Here we insert the intermediate cusps of the chains into the corresponding sides of $F_r$. Note that
 $$
  N_{e_2}=\{5,9,  4\}  \cup \{ 4, 7\},\,\,
 N_{e_3}=\{ 7,9,2\} \cup   \{6,5\},\,\,
 N_{e_4}=\{7,11,4\} \cup \{4,5\},
 N_{e_5}=\{5,11,6\}   \cup  \{2,7\},
 $$
 where $(9,4)$ pairs with $(4,7)$, $(9,2)$ pairs with $(6,5)$, $(11,4)$ pairs with $(4,5)$,
 $(11,6)$ pairs with $(2,7)$.
 In short, \begin{equation}\label{eqn:Sonefor64} S_1=\{ 1,8,7,13,19,6,5,9,4,7,3,8,5,7,9,2,   7,5,8,3,7,11,4,5,11,6,13,7,8,1\}.\end{equation}

\medskip
\noindent 
 In general, $S_1$ is the union of $F_r$ and all the $N_e$'s (intermediate cusps of the $N_e$'s are inserted into the corresponding sides of $F_r$),
  where $e$ is a free side which takes the form $(a, bp^m)$, $m\ge 1$. See Section \ref{s:Sone} for details.

\subsection{(Non)-freeness of sides of $S_1$.}\label{ss:sidesofSone} Notice that all sides of $S_1$ of the form $(a,b)$ where $ab$ is a multiple of $2$ are not free. Put another way, $2$ does not divide $ab$ for any free side $(a,b)$ of $S_1$. This follows from subsection \ref{ss:observation} and the geometric interpretation of subsection \ref{ss:chains}.
We now look at the side $(7,13)$ in $N_{e_1}$ given in (4.4).  This is not a free side.
 We shall in the following give two proofs of this fact. A bad one and a good one that provides some insight.

\smallskip
  \noindent {(B)  \bf A bad proof.}  $(7,13) $ is not free since it pairs with $(7,11) \in S_1$.

  \smallskip
   
   \noindent {(G)  \bf A good proof.}   Suppose that $(7, 13)$ is  free. We may insert $7+13=20$ into $(7,13)$. However, $20$ corresponds to a cusp of intermediate width $w=2^{6-4}=4$ (Lemma \ref{lem:csps2}(i)) but all cusps of intermediate widths have already been filled in by $S_1$, giving a contradiction.
   

     \medskip

     \noindent
     By the same reasoning, if  $(a, b)$ is a side of $S_1$ such that $a+b=z2^i$, $z$ odd, $i < 6/2$,
      then $(a, b)$  is not a free side.   The  only remaining (free) sides are
           $$(13,19)\in N_{e_1}, \,(7,9))\in N_{e_3}, \,(5,11))\in N_{e_5}.$$
      
      \noindent Note that $13+19, 7+9, 5+11$  are multiples of $\sqrt{ 64}$  and that 
      $(13,32)$ pairs with $(32, 19)$, $(7,16)$ pairs with $(16,9)$, $(5,16)$ pairs with $(16,11)$.
        There are no more free sides after the insertion
         of $13\oplus 19, 7\oplus 9, 5\oplus 11$. Hence
       \begin{equation}\label{eqn:Tfor64}T =  S_1\cup \{   13\oplus 19, 7\oplus 9  , 5\oplus 11\}\end{equation}
       \mbox{ has no free sides}.
      
\noindent 
 See subsection \ref{ss:Fareyseq}  for notation. As a consequence, $T$ is a Farey sequence of $\Gamma_0(64)$.      

\subsection{The general case.}\label{ss:generalcase} Let   $S_1$  be the union of $F_r$ and all the $N_e$'s, 
 where $e$ is a free side of the form $(a, b2^m), m\ge 1$. As observed earlier, the sides of $S_1$ of the form
  $(a, b2^m)$ or  $(x2^p, y)$  $(bx$ odd) are not free sides  and the only possible
   free sides of $S_1$   must take the form $(a, b)$, where both $a$ and $b$ are odd.
   Similar to the $\Gamma_0(64)$ case, if $(a,b)$ is a side such that $a+b = z2^i$,
    where $i  < n/2$    ($z$ odd), then $(a, b)$ is not a free side.
     In the case $a+b = z2^i$, $i\ge n/2$, again similar to the $\Gamma_0(64) $ case,
     $(a, a+b)$ and $(a+b, b)$ are paired. In summary,
        \begin{equation}\label{eqn:Tforgeneralcase}
        T =  S_1\cup \{a\oplus b\,:\,  a+b = z2^i, i \ge n/2\}
   \mbox{  has no free sides}.\end{equation}
   Hence $T$ is a Farey sequence of $\Gamma_0(2^n)$. See subsection \ref{ss:pequaltwo} for details and subsection \ref{ss:UpperboundsT} for bounds.

     
\section{ The Farey sequence $S_1$  }\label{s:Sone}

The main theme of this section is  to show that the  Farey sequence $S_1$  constructed from $F_r$ is a Farey subsequence of $\Sigma(p^n)$ as promised in  subsection \ref{ss:SeqSone}. We consider $N=p^n$, $p$ prime, $n \ge 3$, and $r=\lfloor\sqrt{N}\rfloor$ as this works for general prime $p$.

\subsection{ The chains $N_e$.}\label{ss:chainNe}
Let $ e=  (a, b_0p^m)=(a, b)$ be a free side of $F_r$ where
$a$ and $b_0$ are coprime to $p$ and $m\ge 1$.
Consider the chain
\begin{equation}\label{eqn:chain5one} \{ a, a+b_0p^m, a+ 2b_0p^m, \ldots,
a+ vb_0p^m,\ldots,
b_0p^m\}.\end{equation}

\noindent
Since $w(b_0p^m)=p^{n-2m}$ (Lemma \ref{lem:csps2}(i)), there exists a smallest  $0<v<p^{n-2m}$
such that $(a+ vb_0p^m, b_0p^m)$ is paired with a side $(x, y)$ of $F_r$.
Equivalently,
\begin{equation}\label{eqn:chain5end}(a+ vb_0p^m) x + b_0p^m y = kp^n.\end{equation}

\noindent Recall that $b= b_0p^m$. (5.2) therefore  takes
the following simpler form
\begin{equation}\label{eqn:chain5endsimplified}(a+ vb) x + by = kp^n.\end{equation}

\noindent
Note that $x$ must take the form $x_0p^m$, where $x_0$ is coprime to $p$.
Let $v= v_1+v_2$, where $v_1=v_2$ if $v$ is even and $v_1=v_2+1$ if $v$ is odd.
Equation (\ref{eqn:chain5endsimplified}) can be written as
\begin{equation}\label{eqn:chain5endrewritten}(a+ v_1b) x + b(  y + v_2x) = kp^n.\end{equation}
Equation (\ref{eqn:chain5endrewritten}) serves us better as it gives smaller vertices
$(a+v_1b < a+vb$).  We therefore
use (\ref{eqn:chain5endrewritten}) and replace (\ref{eqn:chain5one}) by
\begin{equation}\label{eqn:optimalchains}N_e =   \{ a, a+b,a+2b, \ldots,  a+ v_1b, b\} \cup   \{x, y+v_2x, \ldots, y+x, y
\},\end{equation}
where $(a+v_1b, b)$ pairs with $( x, y+v_2x)$. This gives:

\smallskip

\begin{Lemma}\label{lem:verticesofchains}
    For a free side $e=(a,b_0p^m)$ of $F_r$ as described above, the interior vertices  in $N_e$ are of the form $ a+k_1b$ or $y+k_2x$, where $k_i \le v_i$,
	$v_1+v_2 \le p^{n-m/2}-1$, and  $v_1=v_2$ if $v_1+v_2$ is even and $v_1=v_2+1$ otherwise.
\end{Lemma}

\begin{prop}\label{prop:Sone}
   The congruence subgroup $\Gamma_0(p^n)$ possesses
	a Farey subsequence $S_1$ such that
	
	\begin{enumerate}
		\item[(i)] $S_1=F_r \cup N_e$ where $N_e$ are  the chains of the free sides $e=(a,b_0p^m)$ given in Lemma \ref{lem:verticesofchains} (this union means the interior vertices of $N_e$ are inserted into $F_r$ by (\ref{eqn:Fareyaddition})).
		
		\vspace{.1cm}
		
		\item[(ii)]
		The
		vertices of all the  free sides $(\sigma , \tau)$ of $S_1$
		are coprime to $p$.

\end{enumerate}
\end{prop}

\smallskip
\noindent {\em Proof.} (i) Let $E=\{e_1,e_2,\ldots,e_k\}$ be the  set of free sides of $F_r$ that take the form $(a, b_0p^m)$ as in Lemma \ref{lem:verticesofchains}. We are free to insert
$N_{e_1}$ into $F_r$. By construction, $L_1 =F_r\cup N_{e_1}$ is  a subsequence of a Farey
sequence of $\Gamma_0(p^n)$, since in adding the chain $N_{e_1}$, we only insert into free sides. The key observation here is that when inserting the chain $N_{e_1}$, apart from the two sides $(a+v_1b,b)$ and $(x,v_2y+x)$ which are paired, the other sides introduced are of the form $(u,v)$ where $p$ does not divide $uv$. These cannot be paired to the remaining free $e_j$'s in $E$, as $e_j$ can only be paired with a side of the form $(xp^m,y)$. So $e_2, \ldots, e_k$ remain free. We then look at the second free side $e_2\in E$
of $F_r$ (which is free in $L_1$) which takes the form $(c, dp^m)$ (possibly different $m$ from that for $e_1$) 
and insert $N_{e_2}$ into $L_1$. By the same observation, the sides of the form $(u,v)$ introduced in the insertion of $N_{e_1}$ does not affect the insertion of $N_{e_2}$. Again by the same observation, $e_3, \ldots, e_k$ remain free in $L_2$.
We continue this construction to obtain the Farey sequence
$S_1 = F_r \cup \bigcup_{e_i \in E} N_{e_i}.$ Then, by construction,  $ S_1$ is a subsequence
of a Farey sequence of $\Gamma_0(p^n)$.

\smallskip
\noindent (ii) By construction,   $S_1$ possesses no more free side of the form
$(a, bp^m)$. Using a similar argument, $S_1$ also contains no free sides of the form  $(xp^m, y)$, as this would mean it would have to contain a free side of the form $(a, bp^m)$ when we try to extend the Farey subsequence to a Farey sequence of $\Gamma_0(p^n)$ by successively inserting into the sides having  $xp^m$ as the left vertex. Hence any free side of $S_1$ has the required form. \qed

\subsection{Geometric meaning of filling in the chains}
We note that from Lemma \ref{lem:csps2}(i) and (ii), if $x$ is equivalent to a cusp of the form $a/bp^i$ with $gcd(ab,p)=1$, and $i \le n/2$, then the projection of $F_r$ to $\Sigma(p^n)$ intersects $C(x)$. The gaps in $C(x)$ not in the image are bounded on the left by free sides of $F_r$ of the form $(a, b_0p^m)$ and on the right by free sides of the form $(x_0p^m,y)$ and inserting the chains $N_e$ corresponds to filling in these gaps by inserting Farey triangles to $F_r$. By Lemma \ref{lem:csps2}(iii), only one vertex of each inserted triangle lies in the collar of a cusp of intermediate width, the other two vertices lie in the collar of $[0/1]$, the unique cusp of maximal width. Hence the filling in of each gap in a cusp of intermediate width does not affect the remaining gaps of these cusps, and the chains can be added independently of each other.
Note that in the projection of $S_1$ to $\Sigma(p^n)$, there are no more gaps in the collars of the cusps equivalent to $a/bp^i$, $i \le n/2$.




\section{ Farey sequence $\{e_i\}$ of $\Gamma_0(p^n)$}\label{s:mainconstruction}

\subsection{The case $p=2$.}\label{ss:pequaltwo}
  Let $S_1$ be  given as in Proposition 5.2 and let $e= (a, b)$ be a free
side of $S_1$. 
Following our construction, both $a$ and $b$ are
coprime to $p=2$ (see (ii) of Proposition 5.2).
Hence  $a+b= 2^iz$ is even ($z$ odd).
Since all the cusps equivalent to $2^m z$ are already filled in for $m \le n/2$ by the projection of $S_1$ , $i  >    n/2$. 
As a consequence,
$(a, a+b)$ and  $( a+b, b)$ are paired with each other.
In conclusion,
\begin{equation}\label{eqn:seqforT}
T=S_1\cup \{  a\oplus b\,:\, a+b=z2^i,\,i> n/2 \}  \mbox{ has no free sides}.
\end{equation}
Hence $T$ is a Farey sequence of $\Gamma_0(2^n)$.

\subsection{The general case} Here it may be insufficient to just insert Farey triangles into each of the free sides of $S_1$ to obtain $T$. In general one needs to insert blocks of Farey triangles to the free sides. See Appendix \ref{appendix} for details.

\section{Upper bounds}\label{s:Upperbounds}

 \subsection{Upper bound for vertices in $N_e$}\label{ss:upperboundsforchains}
Let $(a,b) = (a,b_0 p^m) \in F_r$ be a free side, where $ab_0$ is coprime to $p$ and $1 \le m \le n/2$. Then $m < n/2$, since for $m = n/2$ the corresponding cusp at the endpoint with denominator $b_0 p^m$ has width $1$, and its two incident sides are already paired within $F_r$, so none of them can be free.
Let $X= b$ and $Y= p^{ n-m} -a >0$
$(a \le  p^{n/2}, \,m <n/2)$.
Then  $aX +bY =b_0p^{n}$.
Let gcd$\,(X, Y) =d$. Since $Y$ is coprime to $p$,
$d$ is a divisor of $b_0$. In particular, $b_0/d \in \mathbb N$.
We consider
\begin{equation}\label{eqn:lineardiophantineone}a(X/d) + b(Y/d) = (b_0/d)p^n.\end{equation}
This implies that the free side $(a, b)$ of $F_r$ pairs with $(X/d, Y/d)$. Since
$X/d=b/d< p^{n/2}$  and gcd$\,(X/d, Y/d)$
$=1$, one must have $Y/d>p^{n/2}$
(if not, then $(X/d, Y/d) \in F_r$, contradicts the fact that $(a, b)$ is free).
Since  $X/d =b/d < p^{n/2}$,
there exists a unique $w\in \mathbb N$
such that
\begin{equation}\label{eqn:inequalitiesforw1}
0<Y/d- w (X/d) \le  p^{n/2}  < Y/d-(w-1)(X/d).\end{equation}
This implies that \begin{equation}\label{eqn:matchingwithw}
(X/d, Y/d- wX/d  )\in F_r\,\mbox{ and } \,
(a+wb)(X/d)  +b(Y/d-wX/d)=(b_0/d) p^{n}.\end{equation}
Recall that $w >0$.
By   (\ref{eqn:matchingwithw}) and the fact that  $X =b=b_0p^m$,  one has
\begin{equation}\label{eqn:inequalitiesforw2} a+wb=p^{n-m} -(Y-wX) <p^{n-m}
.\end{equation}


\smallskip
\begin{Lemma}\label{lem:boundsforchains}
    
	Let $(a+v_1b)x+b(y+v_2x)=kp^n$ and $w$ be given as in (\ref{eqn:chain5endrewritten}) and  (\ref{eqn:inequalitiesforw1}).
	 Then
  \begin{equation}\label{eqn:inequalitiesforv1andv2}
  a+v_1b\le a+vb < p^{n-1},\qquad    y+v_2x \le y+vx <p^{n-1}.
	\end{equation}
\end{Lemma}
\noindent 
\noindent {\em Proof. }
Let $U = X/d, V=Y/d-wX/d$ and $\tau = b_0/d$ be given as in (\ref{eqn:matchingwithw}). Then 
$(a+ wb)U+ b V = \tau   p^n$. Since $v$ is minimal so that $(a+vb,b)$ is paired to a side of $F_r$, and $(U,V)\in F_r$, $v \le w$. 
Hence $a+v_1b \le a+vb \le a+wb <p^{n-m}$, the last inequality from (\ref{eqn:inequalitiesforw2}). The proof of the second inequality follows similar lines and is left to the reader.\qed


\subsection{Upper bound for vertices  in $T\setminus S_1$} \label{ss:UpperboundsT}
Let  $E$ be the set of free sides of $S_1$.
 Following our result in subsection \ref{ss:generalcase}, details in Section \ref{appendix},  one can obtain a Farey sequence $T$ by  inserting at most $p-1$ vertices into each of the sides $(s,t)$ of $E$. 
   The largest possible vertex one can get is that one does  $p-1$ insertions at a free side, and the insertions
 follow the Fibonacci pattern ($f_1=f_2=1$, $f_m=f_{m-1}+f_{m-2}$) given as follows:
 \begin{equation}\label{eqn:insertionforT}\{ s,  s+t,\ldots\ldots f_{p-1}s+f_{p}t, \ldots \ldots, s+2t, t\}.\end{equation}
  To get an upper bound for all possible free sides, we first assume that  
  $a+(v_1-1)b = s<t=a+v_1b$ (see (\ref{eqn:optimalchains})).  
One can prove easily that the largest vertex  is 
\begin{equation}\label{eqn:Fibonaccibounds}
f_{p-1}s+f_{p}t = f_{p-1}(a+(v_1-1)b)+f_{p}(a+v_{1}b)  \le  (f_{p+1}/2) p^{n-1} 
+ f_pp^{n/2}.  \end{equation}
Note that (\ref{eqn:Fibonaccibounds}) can be obtained from (\ref{eqn:inequalitiesforv1andv2}) and the fact that $v_1 \le (v+1)/2$ since $2v_1=v$ if $v$ is even and $2v_1 -1=v$ if $v$ is odd.
Replace $(s, t) \in N_e$ by $(u, v)\in F_r$ and repeat our construction. Since $u, v\le p^{n/2}$,  one sees clearly that 
$f_{p-1}u +f_p v\le (f_{p+1}/2) p^{n-1} 
+ f_pp^{n/2}. $ In summary, we have:
\begin{prop}\label{prop:boundsforT}
   $T$ can be chosen so that an upper bound for the vertices of $T$ is \begin{equation}\label{eqn:boundsforT}
   (f_{p+1}/2) p^{n-1} 
+ f_pp^{n/2}.\end{equation}
\end{prop}

\section{ The main result}\label{s:mainresults}

Throughout this section, $p=2$ and $n\ge 2$. 
Recall first that $\Gamma_0(2^n)$ has no elliptic elements if $n\ge 2$.  As a consequence, the Farey sequence of $\Gamma_0(2^n)$
  possesses no elliptic sides.

\smallskip
\begin{theorem}\label{thm:mainthm}(The Main Theorem). 
$\Gamma_0(2^n)$
possesses a Farey sequence $ T= \{e_i \}$ such that
\begin{equation}\label{eqn:boundsinmaintheorem} e_i \le 2^{n-1}.\end{equation}
\end{theorem}
\noindent {\em Proof.} 
Let $T= \{e_i\}$ be given as in (\ref{eqn:seqforT}).
Since the upper bound given in (\ref{eqn:boundsforT}) is clearly larger than the upper bound given in (\ref{eqn:inequalitiesforv1andv2}),
 we may assume that the largest vertices come from $T\setminus S_1$.
 By (\ref{eqn:seqforT}),  the largest vertex takes the form $e_j = 2^iz$. Further, (\ref{eqn:boundsforT}) implies that 
\begin{equation}\label{eqn:boundforz} 2^i z  < 2^{n-1} +2^{n/2}.
\,\,\mbox{ Equivalently,  }
z < 2^{ n-i-1} + 2^{n/2-i}.\end{equation}

\noindent Since $ n/2-i < 0$, one has $z \le 2^{n-i-1}$.
As a consequence, $e_j =2^iz \le  2^{n-1}$. \qed

\subsection{The upper bound $2^{n-1}$ is  optimal}\label{ss:upperboundoptimal}
  From Lemmas \ref{lem:cusps1} and \ref{lem:csps2},
  there exists a unique $e_j\in T$ such that
$e_j=z2^{n-1}$. By Theorem \ref{thm:mainthm}, $z=1$. Hence the upper bound $2^{n-1}$ for $T$ is optimal and is attained at this $e_j$.

\subsection{ Independent generators}\label{ss:independentgenerators}  
The set of all side pairings of $T$ forms a set of
independent generators $\{g_i\}$  of $\Gamma_0(2^n)$.
By Theorem \ref{thm:mainthm}, $e_i \le 2^{n-1}$.  The (2,1)-entry of $g_i$ takes the form
$ax+by$, where  $(a,b)$ and $(x,y) $ are sides of $T$  (see \cite{CLLT}).
Hence
an upper bound of the  (2,1)-entries  of the $g_i$'s is
$2\cdot 2^{2n-2}=2^{2n-1}.$

\subsection{Free sides of $S_1$}\label{ss:freesidesofS1}  The following lemma is useful in obtaining 
 a Farey sequence for $\Gamma_0(2^n)$. 

\begin{Lemma}\label{lem:freesidesofNe}
    Let $n\ge3$ and let $(c,d)$ be a free side
	of $S_1$. Then
	 $c$ and $d$ are adjacent vertices of a chain
	\begin{equation}\label{eqn:Neforp2}
 N_e =   \{ a, a+b,a+2b, \ldots,  a+ v_1b, b\}\cup \{ x, y+v_2x, \ldots, y+x, y
	\},\end{equation}
	
	\noindent where $e= (a, b) =(a, b_02^m)\in F_r$ $(ab_0 \,\, odd)$
	and    $(x, y)=(x_02^m, y) \in F_r$  $(yx_0 \,\, odd)$. Furthermore, $m=1$.
	In particular, if $m \ge 2$,
	all the sides $(c, d) $ coming from $N_e$
	are not free sides of $S_1$.
\end{Lemma}

\smallskip
\noindent {\em Proof.}
We first note that if $(a,b)$ is a side of $S_1$ which is also a side of $F_r$, then $(a,b)$ cannot be free so any free side of $S_1$ is a side of $N_e$ for some $e$.
Suppose that
\begin{equation}\label{eqn:freesidesofNe}
c=  a+ kb_02^m,\,\,\, d=  a+ (k+1)b_02^m.\,
\mbox{  Then }    c+d= 2a+(2k+1) b_02^m.\end{equation}

\noindent By
results in subsection \ref{ss:pequaltwo},  $c+d= 2^iz$, $i >n/2\ge 3/2$. Note that
$a$ is odd. Hence  $m=1$.
The case
$c=  y+(k+1)x_02^m,\,\,\, d= y+ kx_02^m$ can be dealt with similarly.
\qed

\section{Algorithm [DKLT] for $\Gamma_0(2^n)$ }\label{s:algo} The following gives our algorithm
 for determining  $T$.

\smallskip
\begin{enumerate}
\item[(i)] Construct $F_r$, where $r=  \lfloor2^{n/2}\rfloor$.

\vspace{.1cm}

\item[(ii)] Check whether $(a, b2^i)$ is free for each such side ($ab$ odd, $i \ge1$).
Note that  $(a, b2^i)$ is not free only if it is paired with
a side of the form $(x2^i, y)$.

\vspace{.1cm}

\item[(iii)] List the free sides in (ii) as $s_1$, $s_2 ,\ldots ,$
$s_k$. We may choose the ordering so that $s_1, s_2, \ldots
,s_r$, $r<k$ are sides of the forms $(a_i, b_i2)$, where $a_ib_i$  is odd. Construct $N_{s_i}$ (see (\ref{eqn:optimalchains})) one by one.



\vspace{.1cm}

\item[(iv)] List all the sides $(a, b)$ coming from inserting the $N_{s_i} \, (i\le r)$,
where $ab$ are odd. Note that Lemma \ref{lem:freesidesofNe}  indicates that
the possible free sides of $S_1$ must come from $N_{s_i}$, $i\le r$.

\vspace{.1cm}

\item[(v)] Let $s= (a, b)$ be a side in (iv).
Insert $a\oplus b$ into
$(a, b)$ if $a+b$ is a multiple of $2^i$  for some  $i\ge n/2 $.

\vspace{.1cm}

\item[(vi)] Since
there are no free sides
left
after the insertions described in
((vi) (see subsection 6.1), the collection of all the
vertices above gives a Farey sequence for $\Gamma_0(2^n)$.
\end{enumerate}

\smallskip
\begin{theorem}\label{thm:algo}
  Algorithm  [DKLT]
 in  this section gives a Farey sequence $T$  with optimal upper bound $2^{n-1}$ for $\Gamma_0(2^n)$.
\end{theorem} 

\smallskip
\noindent {\em Proof. } Apply Theorem \ref{thm:mainthm} and subsection \ref{ss:upperboundoptimal}. \qed

\section{A simple conjecture}
 


\subsection{}
The algorithm we produce in Section 9 can be modified to get a {\em good}  Farey sequence of $\Gamma_0(p^n)$ for many examples 
 which leads us  to 
conjecture that

\medskip

\begin{Conj}\label{Conj:bound}
     If $n\ge 3$, then $\Gamma_0(p^n)$ has a Farey sequence $\{e_i\}$  such that $
 e_i \le p^{n-1}$.
\end{Conj}

\medskip
\noindent  This conjecture is a special case of our more general conjecture for $\Gamma_0(N)$; see \cite[Conjecture 3.16]{DKLT}. In order to prove the conjecture,  we probably need to find a small but pivotal example.
 We have such an example  when  $p=2$ ($\Gamma_0(64)$ for $\Gamma_0(2^n)$).  Unfortunately,
  we don't have such examples when $p$ is odd.

\medskip

  \medskip
  
  \section{Appendix}\label{appendix}
  
  In this section, we show that one may insert at most $p-1$ vertices to each free side of $S_1$ to reach a Farey sequence of $\Gamma_0(p^n)$. We recall the classification of cusps of $\Gamma_0(p^n)$ as follows:
  
\begin{enumerate}
\item  $c=[1/0]=[1/p^n]$; this cusp has width $1$;

\item $c = [a/p^i]$ for each integer $1 \le a < p^i$ with $\gcd(a,p)=1$, and each integer $1 \le i \le n/2$; this cusp has width $p^{n-2i}$;

\item $c = [a/p^i]$ for each integer $1 \le a < p^{n-i}$ with $\gcd(a,p)=1$, and each integer $n/2 < i < n$; this cusp has width $1$;

\item $c = [0/1] = [a/b]$ for every reduced fraction $a/b$ with $\gcd(b,p)=1$; this cusp has width $p^n$.

\end{enumerate}
\noindent
From the discussion in earlier sections, the projection of $S_1$ to $\Sigma(p^n)$ fills the collars of the cusps of types (1) and (2), is empty on  those of type (3), and partially fills the collars of type (4) ($c=[0/1]$). We will need to carefully analyze the structure of the family of Farey triangles around the cusp $c=[0/1]$. 
\subsection{Triangles around the cusp $[0/1]$}
It will be convenient to lift the family of Farey triangles around the cusp $[0/1]$ to the following family of $p^n$ consecutive Farey triangles with a vertex at $0/1$ on $\mathbb{H}$:
$$\Delta :=\left\{ \left(\frac{0}{1},\frac{1}{0}, \frac{1}{1}\right),\left(\frac{0}{1},\frac{1}{1},\frac{1}{2}\right), \ldots , \left(\frac{0}{1}, \frac{1}{p^n-1}, \frac{1}{p^n}\right)\right\},$$
in which each triple is the set of three vertices of the associated triangle. For simplicity, we write
$$\Delta =\{[0,1],[1,2],\ldots,[p^n-1,p^n]\},$$
and implicitly think of each element $[a,b]$ in $\Delta$ as a Farey triangle with vertices $\frac{0}{1},\frac{1}{a}$ and $\frac{1}{b}$. Due to Lemma \ref{lem:csps2}(i)(iii), two distinct elements $[a,a+1]$ and $[b,b+1]$ in $\Delta$ are $\Gamma_0(p^n)$-equivalent, and we write $[a,a+1]\sim [b,b+1]$, if there exists an element $\gamma \in \Gamma_0(p^n)$ that sends the triangle $\left(\frac{0}{1},\frac{1}{a}, \frac{1}{a+1}\right)$ to the triangle $\left(\frac{0}{1},\frac{1}{b}, \frac{1}{b+1}\right)$, in one of the following ways:
\begin{itemize}
\item $\gamma\left(\frac{0}{1}\right)=\frac{1}{b+1}$, $\gamma\left(\frac{1}{a}\right)=\frac{0}{1},$ $\gamma\left(\frac{1}{a+1}\right)=\frac{1}{b}.$
\item $\gamma\left(\frac{0}{1}\right)=\frac{1}{b}$, $\gamma\left(\frac{1}{a}\right)=\frac{1}{b+1},$ $\gamma\left(\frac{1}{a+1}\right)=\frac{0}{1}.$
\end{itemize}
 We have the following properties.
\begin{Lemma}
Two distinct elements $[a, a+1]$ and $[b,b+1]$ in $\Delta$ are equivalent if and only if $a(b+1)+1$ or $b(a+1)+1$ is divisible by $p^n$. In particular, $[a-1, a] \sim [a,a+1]$ if and only if $a^2 \equiv 0$ or $a^2 \equiv -1$ mod $p^n$.

\end{Lemma}
\begin{proof}
Apply the matrix $\begin{bmatrix} a&-1\\ a(b+1)+1 & -b-1 \end{bmatrix}$ or $\begin{bmatrix} a+1 &-1\\ b(a+1)+1 & -b \end{bmatrix}$.
\end{proof}
\noindent
Denote by $N([a,a+1])$ the number of distinct elements in $\Delta$ that are $\Gamma_0(p^n)$-equivalent to $[a,a+1] \in \Delta$.
\begin{Lemma}
For every $[a,a+1] \in \Delta$, we have $N([a, a+1]) \in \{1,2,3\}$, and
\begin{itemize}
\item $N([a, a+1]) = 3$ if and only if $\gcd(a(a+1),p) = 1$, and $a(a+1)\not\equiv -1 \pmod {p^n}$,
\item $N([a, a+1]) = 2$ if and only if $\gcd(a(a+1),p) = p$,
\item $N([a, a+1]) = 1$ if and only if $a(a+1)\equiv -1 \pmod {p^n}$.
\end{itemize}
\end{Lemma}

\begin{proof} 
The result can be deduced from two simple facts:
\begin{itemize}
\item the cusp $[0/1]$ has width $p^n$,
\item a triangle in $\Delta$ has three vertices, with at least two $\Gamma_0(p^n)$-equivalent to the cusp $[0/1]$.
\end{itemize}
\end{proof}

\subsection{Inserting at most $p-1$ vertices}
In $\Delta$, every triangle in the following two disjoint subsets:
\begin{align*}
  U_L &= \left\{[kp^i-1,kp^i] \,\middle|\, 1 \le i \le  n/2,\ 1 \le k \le p^{n-i},\ \text{and }\gcd(k,p)=1 \text{ if } k \ne p^{n-i} \right\}, \\
  U_R &= \left\{[kp^i,kp^i+1] \,\middle|\, 1 \le i \le  n/2,\ 0 \le k < p^{n-i},\ \text{and }\gcd(k,p)=1 \text{ if } k \ne 0 \right\},
\end{align*}
 is $\Gamma_0(p^n)$-equivalent to a triangle in the ideal polygon with vertex set $S_1$, because each such triangle has at least one vertex representing a cusp of type~(1) or type~(2) in the above classification. More precisely, if $\pi:\mathbb{H}\to\Sigma(p^n)$ denotes the projection, then the images of all triangles in $U_L\cup U_R$ lie in $\pi(S_1)$.

\smallskip

\noindent A collection $\{[x,x+1], \ldots , [x+k-1,x+k]\}$ of $k$ consecutive elements of $\Delta$ is called a \textit{gap of length $k$} in $\Delta$, where $k\geq 1$. From the above discussion, we have the following two lemmas. We abuse notation slightly and also denote by $S_1$ the subset of $\Delta$ consisting of triangles equivalent to some triangle associated to $S_1$.
\begin{Lemma}\label{lem: 1}
$\Delta \setminus S_1$ is a collection of gaps of length at most $2p-2$.
\end{Lemma}
\noindent
A gap $\{[x,x+1], \ldots , [x+k-1,x+k]\}$ of $\Delta \setminus S_1$ is \textit{maximal} if both $[x-1,x]$ and $[x+k,x+k+1]$ are not in $\Delta \setminus S_1$.
\begin{Lemma}\label{lem: 2}
If $\{[x,x+1], \ldots , [x+k-1,x+k]\}$  is a maximal gap of $\Delta \setminus S_1$, then the two sides $(\frac{0}{1},\frac{1}{x})$ and $(\frac{0}{1},\frac{1}{x+k})$ are $\Gamma_0(p^n)$-equivalent to two free sides in $S_1$.
\end{Lemma}
\noindent
Our goal is to fill in all the gaps of $\Delta \setminus S_1$ in a reasonable and systematic way. Denote by $$F(S_1) = \{(a_1,b_1),(a_2,b_2),\ldots,(a_u,b_u)\}$$ the set of free sides of $S_1$.

\smallskip

\noindent \underline{Step 1:} We construct $S_2$ as follows:
\begin{itemize}
\item[1.] Insert $a_1+b_1$ to $(a_1,b_1)$ to obtain 
$$S^1_1 =  \{(a_1,a_1+b_1),(a_1+b_1,b_1),(a_2,b_2),\ldots,(a_u,b_u)\}.$$
\item[2.] Consider the next side $(a_2,b_2)$, there are two possibilities:
\begin{itemize} 
\item[i.] If  $(a_2,b_2)$ is still a free side in $S^1_1$. Insert $a_2+b_2$ to $(a_2,b_2)$ to obtain\\ $S^2_1 =  \{(a_1,a_1+b_1),(a_1+b_1,b_1),(a_2,a_2+b_2),(a_2+b_2,b_2),(a_3,b_3),\ldots,(a_u,b_u)\}.$
\item[ii.] If $(a_2,b_2)$ is not a free side in $S^1_1$. We consider the next side $(a_3,b_3)$ instead of $(a_2,b_2)$.
\end{itemize}
\item[3.] In general, $S^{i+1}_1$ is obtained from $S^i_1$ by inserting $a+b$ into a free side $(a,b)$ of $S^i_1$, with an extra condition that $(a,b) \in S_1$. Note that $(a,b)$ is no longer a side in $S^{i+1}_1$.
\item[4.] When our construction ends, we obtain the set $S_2$. Note that in $S_2$, every free side $(a_i,b_i)$ of $S_1$ is either paired or replaced by the pairs $(a_i,a_i+b_i)$ and $(a_i+b_i,b_i)$.
\end{itemize}
\noindent
From the construction of $S_2$, we have the following Lemma:
\begin{Lemma}\label{lem: 3}
$\Delta \setminus S_2$ is a collection of gaps of length at most $2p-4$.
\end{Lemma}
\begin{proof} 
Let $\{[x,x+1], \ldots , [x+k-1,x+k]\}$ be a maximal gap of $\Delta \setminus S_1$. Due to Lemma \ref{lem: 1}, $1\leq k\leq 2p-2$. Due to Lemma \ref{lem: 2}, $(\frac{0}{1},\frac{1}{x})$ is $\Gamma_0(p^n)$-equivalent to a free side, namely $(a_i,b_i)$, in $S_1$. By the above construction, $(a_i,b_i)$ is no longer a free side in $S_2$; thus, the triangle $[x,x+1]$ is no longer an element in $\Delta \setminus S_2$. Similarly, $[x+k-1,x+k]$ is also no longer an element in $\Delta \setminus S_2$. Therefore, if the length of a maximal gap of $\Delta \setminus S_1$ is greater than 1, then in $\Delta \setminus S_2$, it is reduced by at least 2 in length. The result follows.
\end{proof}
\noindent \underline{Step 2:} After Step 1, if $S_2$ still has free sides, we apply Step 1 for $S_2$ and obtain $S_3$. In general, if $S_i$ still has free sides, we follow Step 1 to construct $S_{i+1}$ from $S_i$. The construction will stop after at most $p-1$ steps due to the following lemma:
\begin{Lemma}
$\Delta \setminus S_i$ is a collection of gaps of length at most $2p-2-2i$.
\end{Lemma}
\begin{proof} 
The proof is similar to the proof of Lemma \ref{lem: 3}.
\end{proof}
\noindent
This completes the proof of the bound for $e_i$ for the Farey sequence $\{e_i\}$ of $\Gamma_0(p^n)$ in Theorem \ref{thm:weakerthm}.

\end{document}